\documentclass[11pt]{article}
\usepackage{amsmath,amssymb,latexsym}
\usepackage{theorem}
\newtheorem{theorem}{Theorem}
\newtheorem{lemma}{Lemma}

\newtheorem{corollary}{Corollary}
\newtheorem{proposition}{Proposition}
{\theorembodyfont{\rmfamily} }
\newcommand{\rank}{\operatorname{rank}}

\title{Asymptotics for a generalization of Hermite polynomials}
\author{Manuel Alfaro $^{*}$, Ana Pe\~{n}a $^{*}$ , M. Luisa
Rezola \thanks{Partially supported by MEC of Spain under
Grant MTM2009-12740-C03-03, FEDER funds (EU), and the DGA project
E-64 (Spain).}
\\ Departamento de Matem\'{a}ticas and IUMA. Univ. de Zaragoza (Spain)
\and Juan Jos\'{e} Moreno--Balc\'{a}zar \thanks{Partially supported by MICINN
of Spain under Grant MTM2008--06689--C02--01 and Junta de Andaluc\'{\i}a
(FQM229 and excellence projects FQM481, P06-FQM-1735).}
 \\ Departamento de Estad\'{\i}stica y Matem\'{a}tica Aplicada \\
 Instituto Carlos I de F\'{\i}sica Te\'orica y Computacional \\ Univ. de Almer\'{\i}a (Spain)}
\date{}

\begin{document}

\maketitle

\begin{abstract}
We consider a generalization of the classical Hermite polynomials
by the addition of terms involving derivatives in the inner
product. This type of generalization has been studied in the
literature from the point of view of the algebraic properties.
Thus, our aim is to study the asymptotics of this sequence of
nonstandard orthogonal polynomials. In fact, we obtain
Mehler--Heine type formulas for these polynomials and, as a
consequence, we prove that there exists an acceleration of the
convergence of the smallest positive zeros of these generalized
Hermite polynomials towards the origin.
\end{abstract}


formulas; zeros; Bessel functions.

\bigskip
Corresponding author: Ana Pe\~{n}a, e-mail: anap@unizar.es

Departamento de Matem\'aticas. Universidad de Zaragoza.

50009-Zaragoza  (Spain)

Fax: (+34) 976761338; Phone: (+34) 976761328

\newpage

\section{Introduction}

The study of the asymptotic behavior of orthogonal polynomials is
one of the central problems in the analytic theory of orthogonal
polynomials. In this paper we deal with the asymptotic
properties of some particular families of polynomials orthogonal
with respect to discrete Sobolev inner products of the form:
\begin{equation*}
(P,Q)=\int_I P(x)Q(x)d\mu+\mathbb{P}(c)^t A \mathbb{Q}(c) \,,
\end{equation*}
where $\mu$ is a positive Borel measure supported on an interval
$I \subseteq \mathbb{R}$, $c \in \mathbb{R}$, $A\in
\mathbb{R}^{(s,s)}$ is a positive semidefinite matrix and for a
polynomial with real coefficients $P$ , $\mathbb{P}(x)$ denotes the
column vector  $( P(x), P'(x), \dots, P^{(s-1)}(x))$, being $\mathbb{P}(x)^t$ its transpose. We denote by
$(Q_n)$ the sequence of its orthogonal polynomials and by $(P_n)$
the corresponding one with respect to $\mu$. For general measures
$\mu$, with bounded support and in the Nevai class $M(0,1)$, some asymptotic results can be seen in \cite{mva} and
\cite{lmva}, among others.

A natural approach to study the asymptotic properties of $Q_n$ is
to compare these polynomials with the standard polynomials $P_n$,
whenever the asymptotic properties of these ones are known. In
\cite{alf-mar-rez-ron95} it was studied in detail the asymptotics
of $Q_n$ when $s=2$ and the measure $\mu$ belongs to the Nevai class
$M(0,1),$ so  the support of $\mu$ is bounded. However, when the
support of the measure $\mu$ is unbounded, we have not a so much
general approach (see \cite{MM}). In this framework, the Laguerre
case was considered in \cite{alv-mb}, where
$d\mu=x^{\alpha}e^{-x}dx\, ,$ with $\alpha >-1,$ and $A\in
\mathbb{R}^{(2,2)}$ is a diagonal matrix. There, the authors study
the asymptotic properties of the corresponding orthogonal
polynomials and of their zeros (see also the survey \cite{MM}
where a clearer notation is used). Therefore, the motivation of
this paper is, on the one hand,  to fulfill a gap in the
literature taking into account the other classical measure of
unbounded support $d\mu=e^{-x^2}dx$, and, on the other hand, to
try to obtain a pattern about the scaled asymptotics of
Mehler--Heine type for this unbounded case.

Let us consider the inner product
\begin{equation} \label{pr-hs}
(P,Q)=\int_R P(x)Q(x)e^{-x^2}dx+\mathbb{P}(0)^t A \mathbb{Q}(0),\quad A\in
\mathbb{R}^{(2,2)}\,,
\end{equation}
where the matrix $A$ is positive semidefinite.

The corresponding orthogonal polynomials are somewhat different
from the modified Hermite polynomials considered in
\cite{alv-arv-mar} since in that work  the entry (2,2) in the
matrix $A$ is zero, that is, the term $P'(0)Q'(0)$ does not appear
in the inner product and asymptotic properties were not studied,
either. A particular case of the inner product considered here was
treated in \cite{rec}, where the entry (2,2) in the matrix $A$ is
the only one different from zero. In that paper, the zeros of the
corresponding orthogonal polynomials were studied from the
numerical point of view and a lower bound for the smallest
positive zero was given.

It is important to observe that along this paper  we will focus
our attention on the Mehler--Heine type asymptotics. It is natural
that the reader wonders: \textit{Why is this type of asymptotics
studied?} In an heuristic way it is possible to guess that the
asymptotic behavior of the new orthogonal polynomials and of the
standard orthogonal polynomials is the same in all the complex
plane except for a neighborhood of the origin. Precisely, the
Mehler--Heine type asymptotics describe in a detailed form the
asymptotic behavior around the origin. In this paper, we will
prove that these two families of orthogonal polynomials have a
different behavior around the origin. In this way, Remark 1
of this paper supports the idea that the type of asymptotics that
deserves to be studied is the Mehler--Heine type asymptotics since
it provides the differences between both sequences of orthogonal
polynomials when the degree of the polynomials tends to infinity.

Next, we describe the structure of the paper. In Section 2, we
introduce some properties of the classical Hermite polynomials and
give expressions of the kernel polynomials and their derivatives
which will be used along this paper. In Section 3, we first obtain
Mehler--Heine type formulas for the polynomials orthogonal with
respect to (\ref{pr-hs}) and we observe that the entries in the
matrix $A$ are connected with the order of the Bessel functions
which appear in these formulas. We also notice that the
nondiagonal case does not add further additional information to
the one obtained in the diagonal case concerning this type of
asymptotics (see Theorem \ref{MH2x2}). Furthermore, for the
diagonal case, we deduce how the presence of the masses gives an
asymptotic behavior of the first positive zero different from the
one of the Hermite polynomials. We also observe that the rank of
the matrix does not play an important role to state the
Mehler--Heine type formulas.

These results are a motivation to find a pattern for this type of
asymptotics in a more general framework. Then, looking for this
pattern, by using a symmetrization process
given in \cite{ammr}, in Section 4 we obtain the asymptotic properties of the
polynomials orthogonal with respect to an inner product like
(\ref{pr-hs}) when $A \in \mathbb{R}^{(4,4)}$ is a positive
semidefinite diagonal matrix. More precisely, we prove that the presence
of all the masses produces an increase in four units in the order
of the Bessel function
 appearing in the corresponding Mehler-Heine
type formulas (see Theorem \ref{MH4x0}). In this case, we get a convergence acceleration to 0 of the
two smallest positive zeros. Moreover, we prove that the positive definiteness of the matrix $A$ is a
 necessary condition to
increase the  convergence acceleration of the zeros (see Remark
5). Finally, a conjecture is
raised when $A \in \mathbb{R}^{(2r,2r)}, r \ge 1 \,,$ is a
positive definite diagonal matrix. This conjecture
can be reformulated in terms of the one stated in
\cite{alv-mb} for the Laguerre--Sobolev type polynomials.

\section{Basic tools}
From here on, we will use the standard notation $H_n$ for the
monic classical Hermite polynomials orthogonal with respect to the
weight $e^{-x^2}$. The polynomials  $H_n$ are symmetric, that is,
$H_n(-x)=(-1)^n\,H_n(x)$, and they satisfy (see \cite{sz}):
$$ {\Vert H_n \Vert}^2=\int_R H^2_n(x)\,e^{-x^2}dx=\frac{\sqrt{\pi} n!}{2^n}\,,
 \quad  H'_n(x)=nH_{n-1}(x)\,,$$
$$H_{2n+1}(0)=0\ ,\quad H_{2n}(0)=\frac{(-1)^n(2n)!}{2^{2n} n!}\, .$$

\noindent The $n$th kernel for the Hermite polynomials
$K_n(x,y)=\displaystyle \sum _{k=0}^n \frac{H_k(x)H_k(y)}{\Vert
H_k \Vert ^2}$ satisfies the Christoffel-Darboux formula
\begin{equation*} \label{Ch-D}
K_{n}(x,y)=\frac{1}{{\Vert H_n
\Vert}^2}\frac{H_{n+1}(x)H_n(y)-H_{n+1}(y)H_n(x)}{x-y} \,,
\end{equation*}
from which it follows
$$ K_{2n+1}(x,0)=K_{2n}(x,0)=\frac{
(-1)^{n}}{n!\sqrt{\pi}}\,\frac{H_{2n+1}(x)}{x}\, ,$$ and
$$K_{2n+1}(0,0)=K_{2n}(0,0)=\frac{(2n+1)!}{\sqrt{\pi} 2^{2n}n!^2}\, .$$

\noindent As usual, we denote the derivated kernels  by
$K_n^{(i,j)}(x,y)=\displaystyle\frac{\partial^{i+j}}{\partial{x^i}
\partial{y^j}}K_n(x,y)$

\noindent $=\displaystyle \sum _{k=0}^n
\frac{H_k^{(i)}(x)H_k^{(j)}(y)}{\Vert H_k \Vert ^2}\,(i,j \ge
0)$, with the convention $K_n^{(0,0)}(x,y)=K_n(x,y)$.

 Observe that the symmetry of Hermite polynomials yields to
$K_n^{(i,j)}(0,0)=0$ for all $n$, whenever $i+j$ is an odd integer
number.

In the next lemma we show some formulas for the derivated kernels
that we need throughout the paper

\begin{lemma} \label{nucleos}
The derivated kernels of the Hermite polynomials satisfy
\begin{itemize}
\item[(a)]
\begin{equation*}  K_{2n}^{(0,1)}(x,0)= K_{2n-1}^{(0,1)}(x,0)=\frac{(-1)^{n-1}}{\sqrt{\pi}
(n-1)!}\frac{2xH_{2n}(x)+H_{2n-1}(x)}{x^2}\,,
\end{equation*}

\begin{equation*}  K_{2n-1}^{(0,2)}(x,0)=\frac{2(-1)^{n-1}}{\sqrt{\pi} (n-1)!}\frac{2x H_{2n}(x)+(1-2nx^2)H_{2n-1}(x)}{x^3}\,,
\end{equation*}
\end{itemize}

$ K_{2n}^{(0,3)}(x,0)=\displaystyle {\frac{2(-1)^{n}}{\sqrt{\pi} n!}}\frac{(3-6nx^2) H_{2n+1}(x)-(2n+1)(3-2nx^2)xH_{2n}(x)}{x^4}.$
\begin{itemize}
\item[(b)]
\begin{eqnarray*}
K_{2n-1}^{(1,1)}(0,0)&=&K_{2n}^{(1,1)}(0,0)=\frac{(2n+1)!}{3\,\sqrt{\pi}\,
2^{2n-2}n! (n-1)!}\, ,\\
K_{2n-1}^{(0,2)}(0,0)&=& \frac{-(2n-1)!}{3 \sqrt{\pi}\,2^{2n-4}\,(n-1)!\,(n-2)!} \, ,\\
K_{2n-1}^{(2,2)}(0,0)&=& \frac{(2n-1)!(3n-1)}{15 \sqrt{\pi}\,2^{2n-6}\,(n-1)!\,(n-2)!}\, ,\\
K_{2n}^{(1,3)}(0,0)&=& \frac{-(2n+1)!}{5 \sqrt{\pi}\,2^{2n-4}\,n!\,(n-2)!}\, ,\\
K_{2n}^{(3,3)}(0,0)&=& \frac{(2n+1)!(5n-3)}{35
\sqrt{\pi}\,2^{2n-6}\,n!\,(n-2)!}\, .
\end{eqnarray*}
\end{itemize}
\end{lemma}

\textbf{Proof.} (a) From the Christoffel-Darboux formula taking
successive derivatives with respect to $y$, evaluating at $y=0$,
and using Leibniz's formula, we find for $j=0,1, \cdots$
$$K_n^{(0,j)}(x,0)=\frac{j!}{{\Vert H_n \Vert}^2}\frac{1}{x^{j+1}}\left[
P_j(x,0;H_n)H_{n+1}(x)-P_j(x,0;H_{n+1})H_n(x) \right],$$ where
$P_j(x,0;f)$ is the $\it{j}$th Taylor polynomial of $f$ at $0$ (see, for instance, \cite{alf-mar-rez-ron95}).
Then, the result follows.

(b) To get these formulas it is enough to take into account that,
if for each fixed $j$, we denote by $\sum_{k=j+1}^{n+j+1} \alpha_k
x^k$ the $\it{(n+j+1)}$th Taylor polynomial of $x^{j+1}
K_n^{(0,j)}(x,0)$ at $0$, then $K_n^{(i,j)}(0,0) = i!
\,\alpha_{i+j+1}$. Therefore, from the expressions of the kernels
obtained in (a) and Taylor's formula for the Hermite polynomials,
the result follows after some suitable computations. $\quad \Box$

In the next section we will use the Mehler--Heine type formulas for the monic Hermite
polynomials $H_n$: for $j \in
\mathbb{Z}$ fixed, it holds
\begin{equation}\label{M-H-H2n}
\lim_{n\to \infty}\frac{(-1)^n\,
\sqrt{n}}{n!}H_{2n}\left(\frac{x}{2\sqrt{n+j}} \right)=\left(
\frac{x}{2} \right)^{1/2}\, J_{-1/2}(x)
\end{equation}
\begin{equation}\label{M-H-H2n+1}
\lim_{n\to
\infty}\frac{(-1)^n}{n!}H_{2n+1}\left(\frac{x}{2\sqrt{n+j}}
\right)=\left( \frac{x}{2} \right)^{1/2}\, J_{1/2}(x)
\end{equation}
uniformly on compact sets of the complex plane (see, for instance, \cite[Formulas 22.15.3 and 22.15.4]{ab-st}).

Remind that the Bessel function $J_{\alpha}$ of the first kind of order
$\alpha$ ($\alpha \in \mathbb{R}$) is defined by
\begin{equation*}
J_{\alpha}(z)=\sum_{n=0}^{\infty} \frac{(-1)^n}{n!
\,\Gamma(n+\alpha +1)} \left( \frac{z}{2} \right)^{2n+\alpha},
\end{equation*}
(if $\alpha$ is a negative integer, we assume $n\ge-\alpha$).
Therefore $z^{- \alpha}J_{\alpha}(z)$ is an entire function which
does not vanish at the point $0$.

It is well known that the Bessel functions satisfy the following
recurrence relation (see, for instance, \cite{sz})
\begin{equation}\label{Bessel}
J_{\alpha -1}(z)+J_{\alpha +1}(z)=\frac{2\alpha}{z}\,J_{\alpha}(z)
\end{equation}

In the sequel, the notation $\alpha_n \sim \beta_n$ means that  ${\alpha_n}/{\beta_n}\to 1$ when $n$ goes to
infinity.
\section{Hermite--Sobolev type polynomials}
We denote by $Q_n^{\lambda}$ the monic polynomials orthogonal with
respect to the inner product
\begin{equation} \label{pr-hs1} (P,Q ) =\int_{\mathbb{R}}
P(x)Q(x)e^{-x^2}dx+\mathbb{P}(0)^t A \mathbb{Q}(0)\, ,
\end{equation}
where $A=\left(%
\begin{array}{cc}
  M_0 & \lambda \\
  \lambda & M_1 \\
\end{array}%
\right)$ with $M_0 \ge 0, M_1 \ge 0, \lambda \in \mathbb{R}$  and
$M_0M_1-{\lambda}^2 \ge 0$.

It is obvious that the Hermite polynomials are a particular case
of the polynomials $Q_n^{\lambda}$, as it is enough to take the
matrix $A$ as the zero matrix.

The algebraic properties of these polynomials $Q_n^{\lambda}$ have
been studied in \cite{alf-mar-rez-ron95}. In the next proposition, we
obtain  the estimates for the coefficients which appear in the
representation of the Hermite-Sobolev type  polynomials in terms of
$(H_n)$.

\medskip
We want to make a comment about notation. Obviously, the
polynomials orthogonal with respect to (\ref{pr-hs1}) depend on
the parameters $M_0, M_1$, and $\lambda$, but  there is usually
no confusion denoting them by $Q_n$. However, we want to consider
two different situations, the nondiagonal case $\lambda \not=0$
(so $M_0>0, M_1>0$) and the diagonal case $\lambda =0$ (so $M_0
\ge 0, M_1 \ge 0$). Thus, we do not write the index $\lambda$
whenever  $\lambda =0$, that is, in the diagonal case.

 \medskip
\begin{proposition} \label{pro-mm} Let $Q_n^{\lambda}$ be the monic polynomials orthogonal with respect to the inner product
(\ref{pr-hs1}). Then, for
$n\ge 1,$ the following formulas hold:
\begin{align}
Q^{\lambda}_{2n}(x)&= H_{2n}(x)-a^{\lambda}_n
\,\frac{H_{2n-1}(x)}{x}- b^{\lambda}_n
\,\frac{2xH_{2n}(x)+H_{2n-1}(x)}{x^2} \,
,\label{par} \\
Q^{\lambda}_{2n+1}(x)&= H_{2n+1}(x)-c^{\lambda}_n
\,\frac{H_{2n+1}(x)}{x} -d^{\lambda}_n
\,\frac{2xH_{2n}(x)+H_{2n-1}(x)}{x^2} \label{impar}
\end{align}

where
 $$ \lim_{n\to \infty}a_n \,=\begin{cases} 0 & \text {if }\, M_0=0, M_1 \ge 0 \\
 -\frac{1}{2} & \text {if }\, M_0>0, M_1 \ge 0\end{cases}
\, ; \, \lim_{n\to \infty}a^{\lambda}_n=\begin{cases} 0 & \text {if }\, M_0M_1- {\lambda}^2 = 0 \\
 -\frac{1}{2} & \text {if }\, M_0M_1-{\lambda}^2 > 0 \end{cases} $$

$$\lim_{n\to \infty}\sqrt{n}\,b_n=\lim_{n\to \infty}\sqrt{n}\,b^{\lambda}_n \,=0 \,; \quad \lim_{n\to \infty}\sqrt{n}\,c_n=\lim_{n\to \infty}\sqrt{n}\,c^{\lambda}_n=0$$

$$ \lim_{n\to \infty}d_n \,=\begin{cases} 0 & \text {if }\, M_0\ge 0, M_1 = 0 \\
 -\frac{3}{4} & \text {if }\, M_0 \ge 0, M_1 > 0\end{cases}
\, ; \,\lim_{n\to \infty} d^{\lambda}_n=-\frac{3}{4}  \quad \text {if }\, M_0M_1 - {\lambda}^2 \ge0 \,.$$

\end{proposition}

\textbf{Proof.}  The algebraic relations (\ref{par}) and
(\ref{impar}) can be deduced in the usual way expanding $Q^{\lambda}_n$ in terms of the orthogonal system $(H_n)$ and taking into account the
expressions of $K_n(x,0)$ and $K_n^{(0,1)}(x,0)$, (see also Proposition 6 in \cite{alf-mar-rez-ron95}). The coefficients
are given by
\begin{align*}a^{\lambda}_n&=\frac{(-1)^{n-1}}{\sqrt{\pi}(n-1)!}\frac{H_{2n}(0)}{\Delta^{\lambda}_{2n}}\left[ M_0+(M_0M_1-{\lambda}^2)K_{2n-1}^{(1,1)}(0,0) \right]\\
b^{\lambda}_n&= \frac{(-1)^{n-1}}{\sqrt{\pi}(n-1)!}\frac{H_{2n}(0)}{\Delta^{\lambda}_{2n}}{\lambda} \\
c^{\lambda}_n&=\frac{(-1)^{n}}{\sqrt{\pi}n!}\frac{(2n+1)H_{2n}(0)}{\Delta^{\lambda}_{2n+1}} {\lambda} \\
d^{\lambda}_n
\,&=\frac{(-1)^{n-1}}{\sqrt{\pi}(n-1)!}\frac{(2n+1)H_{2n}(0)}{\Delta^{\lambda}_{2n+1}}\left[M_1+(M_0M_1-{\lambda}^2)K_{2n}(0,0)
\right]\end{align*} where
$$\Delta^{\lambda}_n= 1+M_0K_{n-1}(0,0)+M_1K_{n-1}^{(1,1)}(0,0)+(M_0M_1-{\lambda}^2)K_{n-1}(0,0)K_{n-1}^{(1,1)}(0,0).$$

From the estimates of $H_{2n}(0),K_{n}(0,0), K_{n}^{(1,1)}(0,0)$
and $\Delta^{\lambda}_n$, according to the different cases, and by
using  Stirling's formula adequately, it can be deduced after
suitable computations:
$$a_n=0 \quad \text {if } M_0=0, M_1\ge 0 \, ; \quad \lim_{n\to \infty}a_n= -1/2 \quad \text {if } M_0>0, M_1\ge 0 $$
$$ \lim_{n\to \infty}n\,a^{\lambda}_n=\frac{-3M_0}{8M_1} \quad \text {if } \,  M_0M_1-{\lambda}^2=0$$
$$\lim_{n\to \infty}a^{\lambda}_n=
 -1/2 \quad  \text {if }\, M_0M_1-{\lambda}^2>0 $$
\bigskip
$$b_n=0 \quad \text {if }\, M_0 \ge 0, M_1\ge 0 $$
$$ \lim_{n\to \infty}n\,b^{\lambda}_n=\frac{-3\lambda}{8M_1} \quad \text {if } \quad  M_0M_1-{\lambda}^2=0$$
$$\lim_{n\to \infty}n^{3/2}b^{\lambda}_n=
 \frac{-3 \pi \lambda}{16(M_0M_1-{\lambda}^2)} \quad  \text {if }\quad M_0M_1-{\lambda}^2>0 $$
\medskip
$$c_n=0 \quad \text {if }\, M_0 \ge 0, M_1\ge 0$$
$$ \lim_{n\to \infty}n\,c^{\lambda}_n=\frac{3\lambda}{4M_1} \quad \text {if } \,  M_0M_1-{\lambda}^2=0$$
$$\lim_{n\to \infty}n^{3/2}c^{\lambda}_n=
 \frac{3 \pi \lambda}{8(M_0M_1-{\lambda}^2)} \quad  \text {if }\, M_0M_1-{\lambda}^2>0 $$
\bigskip
$$d_n=0 \quad \text {if }\, M_0 \ge 0, M_1= 0 \, ; \quad \lim_{n\to \infty} d_n=
 -3/4 \quad  \text {if }\, M_0 \ge 0, M_1>0 $$
$$\lim_{n\to \infty} d^{\lambda}_n \,=
 -3/4 \quad  \text {if }\, M_0M_1-{\lambda}^2 \ge 0, $$

and the result follows. $\quad \Box$

\medskip

\noindent \textbf{Remark 1} Taking into account the relative
asymptotics for the monic Hermite polynomials (which can be
obtained from Perron's formula, see \cite[Sect. 8.22]{sz}), i.e.,
$$
\lim_{n \to \infty} \sqrt{n}\,
\frac{H_{2n}(x)}{H_{2n+1}(x)}=-\textrm{sgn}(Im(x))\, i\, ,
$$
 the scaled asymptotics for monic Hermite polynomials (see
\cite[p.126]{wal}), i.e., for $j\in \mathbb{Z}$ fixed and being
$\varphi(x)=x+\sqrt{x^2-1}$ the conformal mapping of
$\mathbb{C}\setminus [-1,1]$ onto the exterior of the closed unit
disk, we have
$$
\lim_{n \to \infty} \sqrt{n}\,
\frac{H_{n-1}(\sqrt{n+j}\,x)}{H_{n}(\sqrt{n+j}\,x)}=\frac{\sqrt{2}}{\varphi\left(x/\sqrt{2}\right)}\,
,
$$ and applying Proposition \ref {pro-mm}, we can deduce after several computations that
$$Q_n^{\lambda}(x)=H_n(x)(1+o(1)),$$
$$ Q_n^{\lambda}(\sqrt{n}\,x)=H_n(\sqrt{n}\,x)(1+o(1)),$$
hold uniformly on compact sets of  $\mathbb{C} \setminus
\mathbb{R}$ and $\mathbb{C} \setminus [-\sqrt 2, \sqrt 2]$,
respectively. Therefore, the polynomials $Q_n^{\lambda}$
have the same outer strong asymptotics and Plancherel--Rotach
type behavior as the Hermite polynomials.

\medskip

We focus our attention on Mehler--Heine type formulas for the polynomials $Q_n^{\lambda}.$ These formulas are interesting
twofold: on the one hand, they provide the scaled asymptotics of
$Q^{\lambda}_n$ on compact sets of the complex plane and, on the
other hand, they supply us with asymptotic information about the location of the zeros of $Q_n^{\lambda}$ in terms of the zeros
of other simple and known special functions.

With the previous results we are ready to prove the scaled
asymptotics for the polynomials $Q^{\lambda}_n$ orthogonal with respect to (\ref{pr-hs1}) where $A$ is not the zero matrix.

\begin{theorem} \label{MH2x2}
The following Mehler-Heine type formulas hold:
\begin{itemize}

\item[(a)] For generalized Hermite polynomials of even  degree

$$\lim_{n\to \infty}\frac{(-1)^n\, \sqrt{n}}{n!}Q_{2n}\left(\frac{x}{2\sqrt{n}}
\right)=\begin{cases} (\frac{x}{2})^{1/2} J_{-1/2}(x) & \text {if $\rank A=1, M_0=0$}\,, \\
-(\frac{x}{2})^{1/2} J_{3/2}(x) & \text {if $\rank A=1, M_0>0$}\,, \\
 -(\frac{x}{2})^{1/2} J_{3/2}(x) & \text {if $\rank A=2$ } \, . \end{cases}$$

$$\lim_{n\to \infty}\frac{(-1)^n\, \sqrt{n}}{n!}Q_{2n}^{\lambda}\left(\frac{x}{2\sqrt{n}}
\right)=\begin{cases}  (\frac{x}{2})^{1/2} J_{-1/2}(x) & \text {if $\rank A=1$} \, ,\\
-(\frac{x}{2})^{1/2} J_{3/2}(x) & \text {if $\rank A=2$} \, .\end{cases}$$

\item [(b)]For generalized Hermite polynomials of odd degree

$$\lim_{n\to \infty}\frac{(-1)^n}{n!}Q_{2n+1}\left(\frac{x}{2\sqrt{n}}
\right)=\begin{cases} (\frac{x}{2})^{1/2} J_{1/2}(x) & \text {if $\rank A=1, M_1=0$} \, ,\\
-(\frac{x}{2})^{1/2} J_{5/2}(x) & \text {if $\rank A=1, M_1>0$}\, ,\\
-(\frac{x}{2})^{1/2} J_{5/2}(x) & \text {if $\rank A=2$ }
\, .\end{cases}$$

$$\lim_{n\to \infty}\frac{(-1)^n}{n!}Q_{2n+1}^{\lambda}\left(\frac{x}{2\sqrt{n}}
\right)=\begin{cases} -(\frac{x}{2})^{1/2} J_{5/2}(x) & \text {if $\rank A=1$} \, ,\\
-(\frac{x}{2})^{1/2} J_{5/2}(x) & \text {if $\rank A=2$} \, .\end{cases}$$

\end{itemize}

All the limits hold uniformly on compact sets of the complex
plane.
\end{theorem}

\textbf{Proof.} (a) We have only written the proof for the
nondiagonal case $(\lambda \not= 0)$, since the diagonal case
$(\lambda = 0)$ can be deduced in a similar way. From formula
(\ref{par})  in Proposition \ref{pro-mm}, we have
\begin{align*}&\frac{(-1)^n\,\sqrt{n}}{n!}Q^{\lambda}_{2n}\left(\frac{x}{2\sqrt{n}}
\right)=\frac{(-1)^n\sqrt{n}}{n!}H_{2n}\left(\frac{x}{2\sqrt{n}}
\right)\\&+ 2a^{\lambda}_n\,\frac{1}{x}
\frac{(-1)^{n-1}}{(n-1)!}H_{2n-1}\left(\frac{x}{2\sqrt{n}}
\right)\\&+4b^{\lambda}_n\sqrt{n}\left[
\frac{1}{x^2}\frac{(-1)^{n-1}}{(n-1)!}H_{2n-1}\left(\frac{x}{2\sqrt{n}}\right)-\frac{1}{x}\frac{(-1)^{n}\sqrt{n}}{n!}H_{2n}\left(\frac{x}{2\sqrt{n}}
\right) \right]\, . \end{align*}

Now, we must distinguish two different cases according to the
asymptotic behavior of the coefficients $a^{\lambda}_n$ and
$b^{\lambda}_n$. Handling  formulas  (\ref{M-H-H2n}),
(\ref{M-H-H2n+1}), and (\ref{Bessel}) adequately, we get:

\begin{itemize}
\item {} If  $\rank A=1$ (that is $M_0M_1-\lambda^2=0, M_0>0, M_1>0$), since $a^{\lambda}_n \to 0$ and
$\sqrt{n}\,b^{\lambda}_n \to 0$, then
$$\lim_{n\to \infty}\frac{(-1)^n\sqrt{n}}{n!}Q^{\lambda}_{2n}\left(\frac{x}{2\sqrt{n}}
\right) =\left( \frac{x}{2} \right)^{1/2}\, J_{-1/2}(x) .$$ \item
{} If $\rank A=2$ (that is $M_0M_1-\lambda^2>0, M_0>0, M_1>0$), since $2a^{\lambda}_n \to -1$ and
$\sqrt{n}\,b^{\lambda}_n \to 0$, then
\begin{align*}&\lim_{n\to \infty}\frac{(-1)^n\sqrt{n}}{n!}Q^{\lambda}_{2n}\left(\frac{x}{2\sqrt{n}}
\right) \\&= \left( \frac{x}{2} \right)^{1/2}\left[J_{-1/2}(x)
-\frac{1}{x}J_{1/2}(x)\right]=-\left( \frac{x}{2} \right)^{1/2}\,
J_{3/2}(x).\end{align*} All the limits hold uniformly on compact
sets of $\mathbb{C}$.
\end{itemize}

(b) Upon in the case $M_1=0$ and $\lambda = 0$
($Q_{2n+1}=H_{2n+1}$), the asymptotic behavior of
$\sqrt{n}\,c^{\lambda}_n$ and $d^{\lambda}_n \,$ is independent
from the occurrence of the parameter $\lambda$ in the inner
product, where the coefficients $c^{\lambda}_n$ and $d^{\lambda}_n
\,$ are those ones appearing in formula (\ref{impar}). So, we do
not write the index $\lambda$. Since $\sqrt{n}\,c_n \to 0$ and
$4\,d_n \to -3$, using (\ref{M-H-H2n}), (\ref{M-H-H2n+1}) and
property (\ref{Bessel}), we get:
\begin{align*}&\lim_{n\to \infty}\frac{(-1)^n}{n!}Q_{2n+1}\left(\frac{x}{2\sqrt{n}}
\right)=\lim_{n\to \infty}\left\{
\frac{(-1)^n}{n!}H_{2n+1}\left(\frac{x}{2\sqrt{n}} \right)
\right.\\&-2\sqrt{n}\,c_n\,\frac{1}{x}
\frac{(-1)^{n}}{n!}H_{2n+1}\left(\frac{x}{2\sqrt{n}} \right)\\&
\left. -4\,d_n \left[ \frac{1}{x}\frac{(-1)^{n}\,
\sqrt{n}}{n!}H_{2n}\left(\frac{x}{2\sqrt{n}}
\right)-\frac{1}{x^2}\frac{(-1)^{n-1}}{(n-1)!}H_{2n-1}\left(\frac{x}{2\sqrt{n}}
\right) \right] \right\} \\
&=\left(\frac{x}{2}\right)^{1/2} \left[ J_{1/2}(x)+\frac{3}{x} \left( J_{-1/2}(x) -\frac{1}{x}J_{1/2}(x)\right) \right]\\
&=\left(\frac{x}{2}\right)^{1/2} \left[ J_{1/2}(x)-\frac{3}{x}
J_{3/2}(x) \right]=-\left(\frac{x}{2}\right)^{1/2}J_{5/2}(x)\,,
\end{align*}
holds uniformly on compact sets of $\mathbb{C}$. $\quad \Box$

\medskip

\noindent \textbf{Remark 2} Observe that the nondiagonal case does
not add further additional information to the one obtained in the
diagonal case concerning the asymptotic behavior of the scaled
polynomials.

For this reason, from now on, we will only consider the inner
product (\ref{pr-hs1}) with $\lambda=0$, i.e., the diagonal case.

 \medskip

\noindent \textbf{Remark 3} We want to emphasize the fact that
what it is really important in order
  to have a different result from the Hermite case is the
presence of the masses either $M_0$ for the polynomials $Q_{2n}$
or $M_1$ for  $Q_{2n+1}$ and not the $\rank$ of $A$. Observe how the presence of these masses
implies, in addition to a change of sign, an increase in the order
of the corresponding Bessel functions appearing in Theorem
\ref{MH2x2}.

\medskip

On the other hand, in the next corollary we will show a remarkable difference between the zeros of $(H_n)$ and
$(Q_n)$ with respect to the convergence acceleration to $0$.

 Before analyzing this, we recall (see \cite{sz}) that the zeros of the Hermite polynomials are real, simple
and symmetric. We denote by $(x_{n,k})_{k=1}^{[n/2]}$ the positive
ones in increasing order. It is worth pointing out that they satisfy the
interlacing property $0 < x_{n+1,1} < x_{n,1} < x_{n+1,2} <
\dots $, and that $x_{n,k} \underset{n} \to 0$ for
every fixed $k$.

 Let $(j_{\alpha,k})_{k\ge 1}$ be the positive zeros of the Bessel function $J_{\alpha}$ in increasing order.
Then, formulas (\ref{M-H-H2n}) and (\ref{M-H-H2n+1}) and Hurwitz's
theorem lead us to
 $$2 \sqrt{n}x_{2n,k}\underset{n} \to j_{-1/2,k} \quad (k\ge 1)$$
 $$2 \sqrt{n}x_{2n+1,k} \underset{n} \to j_{1/2,k} \quad (k\ge 1)$$
 and therefore
$$x_{n,k} \sim \frac{C_k}{\sqrt n} \quad (k\ge 1)$$
where, for every $k$, $C_k$ is a positive constant.

Concerning the zeros of $Q_n$ we know  that all of them are real,
simple and symmetric and they interlace with those of $H_n$
(see \cite{alf-mar-rez-ron92}). We denote by
$(\xi_{n,k})_{k=1}^{[n/2]}$ the positive zeros of  $Q_n$ in
increasing order. In this case, it also happens that $\xi_{n,k} \underset{n} \to 0$ for each fixed $k$.

From Theorem \ref{MH2x2} and Hurwitz's theorem and taking into account the multiplicity of 0 as a
zero of the limit functions in Theorem \ref{MH2x2} we achieve

\begin{corollary} Let $(\xi_{n,k})_{k=1}^{[n/2]}$ be the positive zeros of  $Q_n$ in increasing order. Then

\begin{itemize}
  \item [(a)] If $M_0=0$
$$2 \sqrt{n}\xi_{2n,k} \underset{n} \to j_{-1/2,k} \quad (k\ge 1).$$

If $M_0>0$
$$\sqrt n \, \xi_{2n,1} \underset{n} \to 0$$
$$2 \sqrt{n}\xi_{2n,k} \underset{n} \to j_{3/2,k-1} \quad (k\ge 2).$$

\item[(b)] If $M_1=0$
$$2 \sqrt{n}\xi_{2n+1,k} \underset{n} \to j_{1/2,k} \quad (k\ge 1).$$

If $M_1> 0$
$$\sqrt n \, \xi_{2n+1,1} \underset{n} \to 0$$
$$2 \sqrt{n}\xi_{2n+1,k} \underset{n} \to j_{5/2,k-1} \quad (k\ge 2).$$

\end{itemize}
\end{corollary}

Observe that in all the cases we have $\displaystyle \xi_{n,k}
\sim \frac{C_k}{\sqrt n} \quad (k\ge 2).$ However, there exist only
two situations for which the asymptotic behavior of the first
positive zero is different from the one of the Hermite
polynomials. They  correspond  to $M_0>0$ for even degree
polynomials and $M_1>0$ for odd degree polynomials, then $\sqrt n
\, \xi_{n,1} \underset{n} \to 0$. Thus, the presence of the masses
$M_0$ and $M_1$ in the inner product (\ref{pr-hs1}) produces a
convergence acceleration to $0$ of two zeros of the polynomials
$(Q_n)$, namely, the first positive zero and its symmetric one.

\section{Mehler--Heine type formulas: The diagonal case}

The comments and results in the previous section are a motivation
to study what happens with these properties when the matrix$A$
is a positive semidefinite and diagonal .

We begin considering a  diagonal matrix  $A\in
\mathbb{R}^{(4,4)}$. Thus, we introduce the inner product
\begin{equation} \label{pr-hs2} (P,Q) =\int_{\mathbb{R}}
P(x)Q(x)e^{-x^2}dx+\sum_{i=0}^3 M_i\,P^{(i)}(0)\,Q^{(i)}(0) ,
\end{equation}with $M_i \ge 0, \quad i=0,1,2,3$.
We denote by $S_n$ the monic orthogonal polynomials with respect
to (\ref{pr-hs2}).

Notice that in this case the polynomials $S_n$ are symmetric,
i.e., $S_n(-x)=(-1)^n\,S_n(x).$ This doesn't occur for the polynomials
$Q_n^{\lambda}$ considered in the previous section when
$\lambda\neq 0.$ Therefore, because of this symmetry, we can transform the inner product
(\ref{pr-hs2}) into a Laguerre--Sobolev type inner product and so
we can establish a simple relation between the polynomials $S_n$
and the polynomials studied in \cite{alv-mb} and \cite{MM}. This
technique is known as a symmetrization process. In fact, in \cite{ch} this
process is considered for standard inner products associated with
positive measures. The simplest case of this
situation is the relation between  monic Laguerre polynomials and
Hermite polynomials, that is (see \cite{ch} or \cite{sz}),

$$
H_{2n}(x)=L_n^{(-1/2)}(x^2), \quad H_{2n+1}(x)=xL_n^{(1/2)}(x^2),
\quad n\ge 0.
$$
 Later in \cite{ammr} the authors generalize the symmetrization process in the framework of Sobolev type
 orthogonal polynomials.

 Thus, applying Theorem 2 in \cite{ammr} in a straightforward way
 we obtain that $$S_{2n}(x)=L_n^{(-1/2, M_0,4M_2)}(x^2) \quad \mathrm{and} \quad S_{2n+1}(x)
 =xL_n^{(1/2, M_1,36M_3)}(x^2), $$
where $\left(L_n^{(-1/2, M_0,4M_2)}\right)$ and $\left(L_n^{(1/2,
M_1,36M_3)}\right)$ are the sequences of monic orthogonal
polynomials with respect to
\begin{align*}
(P,Q)_{1}&=\int_{0}^{\infty} P(x)Q(x)
x^{-1/2}e^{-x}dx+M_0P(0)Q(0)+4M_2P'(0)Q'(0),\\
(P,Q)_{1^*}&=\int_{0}^{\infty} P(x)Q(x)
x^{1/2}e^{-x}dx+M_1P(0)Q(0)+36M_3P'(0)Q'(0),
\end{align*}
respectively. The Mehler--Heine type formulas for the orthogonal
polynomials with respect to the above inner products were obtained
in \cite{alv-mb} and later reformulated more clearly in \cite{MM}.
Observe that the inner products considered in those articles  are
$\displaystyle{\frac{(P,Q)_{1}}{\sqrt{\pi}}}$ and
$\displaystyle{\frac{2(P,Q)_{1^*}}{\sqrt{\pi}}}$, respectively.
Taking into account that the Mehler--Heine type formulas do not
depend on the explicit value of the masses $M_i$, but  only on
whether the masses are positive or not, see Proposition 2.10 in
\cite{MM}, we deduce the following result directly:

\begin{theorem}  \label{MH4x0}The polynomials $S_n$ satisfy the following Mehler-Heine type formulas.
\begin{itemize}
\item[(a)] For polynomials of even  degree:

\begin{align*}
&\lim_{n\to \infty}\frac{(-1)^n\,
\sqrt{n}}{n!}S_{2n}\left(\frac{x}{2\sqrt{n}}
\right)\\&=\begin{cases} (\frac{x}{2})^{1/2} J_{-1/2}(x) & \text {if }\, M_0=0, M_2=0\, , \\
-(\frac{x}{2})^{1/2} J_{3/2}(x) & \text {if }\, M_0>0, M_2=0 \, ,\\
 (\frac{x}{2})^{1/2} \left[\frac{2}{3}J_{7/2}(x)- J_{3/2}(x)-\frac{2}{3}J_{-1/2}(x)\right] & \text {if }\, M_0=0, M_2>0 \, ,\\
 (\frac{x}{2})^{1/2} J_{7/2}(x)& \text {if }\, M_0>0, M_2 >0\, .\\
 \end{cases} \end{align*}

\item [(b)]For polynomials of odd degree:

\begin{align*}
&\lim_{n\to
\infty}\frac{(-1)^n}{n!}S_{2n+1}\left(\frac{x}{2\sqrt{n}}
\right)\\&=\begin{cases} (\frac{x}{2})^{1/2} J_{1/2}(x) & \text {if }\, M_1=0, M_3=0 \, ,\\
-(\frac{x}{2})^{1/2} J_{5/2}(x) & \text {if }\, M_1>0, M_3=0 \, ,\\
 (\frac{x}{2})^{1/2} \left[\frac{2}{5}J_{9/2}(x)- J_{5/2}(x)-\frac{2}{5}J_{1/2}(x)\right] & \text {if }\, M_1=0, M_3>0 \, ,\\
 (\frac{x}{2})^{1/2} J_{9/2}(x)& \text {if }\, M_1>0, M_3 >0\, .\\
 \end{cases} \end{align*}

\end{itemize}
All the limits hold uniformly on compact sets of the complex
plane.
\end{theorem}

\medskip

Notice that in the even case as well as in the odd case, the
presence of the two relevant masses produces an increase in four
units in the order of the Bessel functions appearing in the
corresponding Mehler--Heine type formulas. As we have proved, we
get an increase in two units in the order of Bessel functions,
when only the first masses appear in (\ref{pr-hs2}), i.e., $M_0>0$
and $M_2=0$ or $M_1>0$ and $M_3=0$ in the respective cases.

\medskip

\noindent \textbf{Remark 4} We can also observe that the rank of
the matrix $A$ is not relevant to establish the Mehler--Heine type
formulas as we have said in the previous section.

In the next corollary we only display the results which are
different from those ones obtained before.

\begin{corollary} Let $(\xi_{n,k})_{k=1}^{[n/2]}$ be the positive zeros of $S_n$ in increasing order. Then
\begin{itemize}

\item[(a)] If $M_0=0$ and $M_2 > 0, \quad  \xi_{2n,k} \sim
\displaystyle \frac{C_k}{\sqrt{n}} \quad (k\ge 1)$

If $M_0>0$ and $M_2 > 0$
$$\sqrt n \, \xi_{2n,k} \underset{n} \to 0 \quad (k = 1, 2)$$
$$\xi_{2n,k} \sim \frac{C_k}{\sqrt{n}} \quad (k\ge 3)$$

\item[(b)] If $M_1 = 0$ and $M_3>0, \quad \xi_{2n+1,k} \sim
\displaystyle \frac{C_k}{\sqrt{n}} \quad (k\ge 1)$

If $M_1> 0$ and $M_3> 0$
$$\sqrt n \, \xi_{2n+1,k} \underset{n} \to 0 \quad (k = 1, 2)$$
$$ \xi_{2n+1,k} \sim \frac{C_k}{\sqrt{n}} \quad (k\ge 3) \,.$$
\end{itemize}
In all the cases, for every $k$, $C_k$ is a positive constant.
\end{corollary}

\bigskip
\noindent \textbf{Remark 5} We want to point out that a singular
fact occurs when there is a gap in the set of the masses, namely
$M_0 = 0\,, M_2 > 0$ or $M_1 = 0\,, M_3 > 0$, in the respective
cases. This difference appears as much in the expression of the
limit function (a particular linear combination of Bessel
functions, see Theorem \ref{MH4x0}) as in the convergence
acceleration to $0$ of the zeros. Observe that in order to get a
convergence acceleration to $0$ of four zeros of $(S_n)$ (the two
smallest positive zeros and their symmetric ones) it is necessary
that all the masses $M_i$
 appearing in the inner product (\ref{pr-hs2}) are positive.

\bigskip

In \cite{alv-mb} a nice conjecture was stated for the orthogonal
polynomials with respect to a Laguerre--Sobolev type inner product
involving $r$ masses at the origin. This conjecture was
reformulated with a clearer notation in the survey paper
\cite{MM}. Therefore, according to our previous results it is
natural to pose a similar one for the orthogonal polynomials,
$Q_n,$ with respect to the inner product

\begin{equation*}
 (P,Q) =\int_{\mathbb{R}}
P(x)Q(x)e^{-x^2}dx+\sum_{i=0}^{2r-1} M_i\,P^{(i)}(0)\,Q^{(i)}(0) ,
\quad r \ge 1, \, M_i \ge 0 \, .
\end{equation*}

Then, it should be true:

\medskip
\noindent \textbf{Conjecture.} If $M_i >0 \,, i=0,1,\dots, 2r-1,$
with $r \ge 1$, then

\begin{equation*} \label{conj1}
\lim_{n\to \infty}\frac{(-1)^n\,
\sqrt{n}}{n!}Q_{2n}\left(\frac{x}{2\sqrt{n}} \right)= (-1)^{r}
\left( \frac{x}{2} \right)^{1/2}\, J_{-\frac{1}{2} +2r}(x)
\end{equation*}
\begin{equation*}\label{conj2}
\lim_{n\to
\infty}\frac{(-1)^n}{n!}Q_{2n+1}\left(\frac{x}{2\sqrt{n}} \right)=
(-1)^{r} \left( \frac{x}{2} \right)^{1/2}\, J_{\frac{1}{2}+2r}(x)
\,.
\end{equation*}
uniformly on compact sets of $\mathbb{C}$.

\medskip

 Now, by using again the symmetrization process given in \cite[Theorem
 2]{ammr}, we can rewrite  our polynomials $Q_n$ as

$$ Q_{2n}(x)=L_n^{(-1/2, N_0,\ldots, N_{2r-2})}(x^2) \quad \mathrm{and} \quad Q_{2n+1}(x)
 =xL_n^{(1/2, N_1,\ldots, N_{2r-1})}(x^2), $$
where $\left(L_n^{(-1/2, N_0,\ldots, N_{2r-2})}\right)$ and
$\left(L_n^{(1/2, N_1,\ldots, N_{2r-1})}\right)$ are the sequences
of monic orthogonal polynomials with respect to
\begin{align*}
(P,Q)_{r-1}&=\int_{0}^{\infty} P(x)Q(x)
x^{-1/2}e^{-x}dx+\sum_{i=0}^{r-1}N_{2i}P^{(i)}(0)Q^{(i)}(0),\\
(P,Q)_{{(r-1)}^*}&=\int_{0}^{\infty} P(x)Q(x)
x^{1/2}e^{-x}dx+\sum_{i=0}^{r-1}N_{2i+1}P^{(i)}(0)Q^{(i)}(0),
\end{align*}
where $$N_0=M_0, \quad N_{2i}=(i+1)_i^2M_{2i}\quad \mathrm{and}
\quad N_{2i+1}=(i+1)_{i+1}^2 M_{2i+1}$$ and $(a)_i$ denotes the
Pochhammer symbol, that is, $$(a)_i=a(a+1)\cdots
(a+i-1)\, ,\quad (a)_0=1.$$

We have proved this conjecture for $r=1$ and $r=2$ in Theorems
\ref{MH2x2} and \ref{MH4x0}, respectively. However, the techniques
used in this paper and in \cite{alv-mb}  do not seem the most
adequate ones to prove the general case. We want to
highlight that in solving the conjecture for the Laguerre case in
\cite{alv-mb} we have solved the one for the Hermite case.

Finally, it is worth observing that as a consequence the following
result could be deduced for the positive zeros of $Q_n$ :

$$\sqrt n \, \xi_{n,k} \underset{n} \to 0 \quad (k = 1, 2, \dots, r)$$
$$\xi_{n,k} \sim \frac{C_k}{\sqrt{n}} \quad (k\ge r+1)\,.$$

So, the presence of all the constants $M_i >0 \,, i=0,1,\dots,
2r-1,$ in the above generalized inner product would induce a
convergence acceleration to $0$ of $2r$ zeros of the polynomials
$(Q_n)$, namely, the $r$ smallest positive zeros and their
symmetric ones.

\medskip

\noindent \textbf{Acknowledgments.} The authors thank one of the referees
for his relevant suggestions which have made
the paper shorter and more readable.

\end{document}